\documentclass[onecolumn,10pt]{article}
\usepackage{algorithm, algorithmic}
\usepackage{multirow,bigdelim}
\usepackage{amsmath, amssymb}
\usepackage{nicematrix}
\usepackage{bm}
\usepackage{hyperref, cleveref}
\usepackage{amscd}
\usepackage{xurl}
\usepackage[numbers,square,sort]{natbib}
\usepackage{doi}
\usepackage{graphicx,subcaption}
\usepackage{xcolor}
\usepackage{amsfonts} 
\usepackage{fullpage} 
\usepackage[tableposition=top]{caption}
\newtheorem{define}{Definition}[section]

\newtheorem{prop}[define]{Proposition}

\makeatletter
   
   \@addtoreset{table}{section}
  \makeatother

\hypersetup{
    breaklinks=true
}

\usepackage{booktabs}
\usepackage{tabularx}
\usepackage{amsmath}
\usepackage{graphicx}

\newcommand{\R}{\mathbb{R}}
\newcommand{\Rm}[2]{\mathbb{R}^{#1 \times #2}}
\newcommand{\avg}[1]{\mathbb{E}\left[#1\right]}
\newcommand{\var}[1]{\mathbb{V}[#1]}
\newcommand{\email}[1]{\href{mailto:#1}{#1}}

\usepackage{tikz}
\usetikzlibrary{arrows.meta,positioning}
\usepackage{nicematrix}

\title{A multilevel sketch-and-solve method for overdetermined least squares problems
} 

\author{Irina-Beatrice Haas,  Michael B. Giles, and Yuji Nakatsukasa \thanks{Mathematical Institute, University of Oxford, OX2 6GG, UK (\email{haas@maths.ox.ac.uk}, \email{mike.giles@maths.ox.ac.uk}, \email{nakatsukasa@maths.ox.ac.uk}).}}
\date{\today}
\begin{document}
\large
\maketitle

\begin{abstract}
    Sketch-and-solve (SAS) is a very successful method to efficiently estimate the solution of heavily overdetermined large linear least squares problems. It uses random sketching to reduce the size of the problem, hence reducing the computational cost.
    Several authors have shown that averaging several solutions from SAS further improves the accuracy, which is measured by the residual associated to the approximate solution.

    Going further, we combine solutions from sketch-and-solve in a multilevel manner, such that the approximate solution is a combination of SAS samples obtained from small sketches and more accurate correction terms obtained from larger sketches. We first consider the variance of the estimator, which depends on the variance of the coarse samples and the correction terms. We show that the variance of the correction terms on each level follows a trend and decreases faster than the variance of the simple SAS estimator. However, we then show that the overall computational cost of our multilevel framework is slightly higher than that of the simple average estimator, so a naive application of multilevel methods appears unattractive for least squares problems. 
     
\end{abstract}



\section{Introduction}

Solving large-scale overdetermined least squares problems is an important step in many optimization tasks and applications. Classically the problem
\begin{equation}\label{eq:LSpb}
    \min_x \|Ax-b\|_2 
\end{equation}
where $A\in \Rm{m}{n}$ and $b\in \R^m$ (for $m>n$), is solved using the QR decomposition of the matrix $A$ with a cost of $\mathcal{O}(mn^2)$ flops. When $m$ is extremely large and $m\gg n$, one can reduce the size of the problem to $s\times n$ by sketching the problem \eqref{eq:LSpb} on the left, (with) the usual choice for the \textit{sketch size} being $s= cn$, where for example $c=4$. Then the solution $\hat x$ of the sketched problem
\begin{equation} \label{eq:SASpb}
    \min_x \|S(Ax-b)\|_2.   
\end{equation}
can be computed in $\mathcal{O}(sn^2)$ operations and gives a good estimate of the exact solution $x^*$ of \eqref{eq:LSpb}. The \textit{sketch matrix} $S\in \Rm{s}{m}$ is chosen in practice such that: (1) it preserves distances up to a factor, which leads to accuracy guarantees for $\hat x$, and (2) the products $SA$ and $Sb$ are cheap to compute. Much of the existing literature has focused on constructing such matrices; see for example \cite{woodruff2014sketching} and the references therein and Section \ref{sec:BackgroundSketching} for an account of the methods that we consider in this paper.

The sketch-and-solve (SAS) paradigm was first introduced 
in the computer science community \cite{drineas2006sampling,drineas2011faster} where it was applied to low-rank approximation of matrices, least squares problems and fast matrix multiplication \cite{sarlos,drineas2011faster}. 
The framework has been analyzed for different sketching methods in terms of the accuracy of the solution $\hat x$ as an approximation of $x^*$ \cite{mahoney2011randomized,woodruff2014sketching,drineasRandNLA}, 
showing that for an appropriate choice of the sketch matrix $S$, with high probability the residual satisfies
\begin{equation}\label{eq:resBound}
    \|A\hat x-b\|_2^2\leq (1+\eta) \|Ax^*-b\|_2^2.
\end{equation}
Here $\eta$ denotes the subspace embedding constant, and usually not very small, e.g. $\eta=0.5$.

In scientific computing, such heavily overdetermined least squares problems arise naturally in the context of subspace methods for linear systems~\cite{nakatsukasa_tropp_2024fast} or model order reduction~\cite{jones2025subapsnap}, where sketch-and-solve has been successfully applied. 

In statistics, several authors \cite{ma15a,raskutti,pilanci2016iterative,wang17d, sketched_ridgeReg} consider the problem of parameter estimation for a data set that follows a Gaussian linear regression model 
\begin{equation}\label{eq:model}
    b = Ax_{pop}+\xi
\end{equation}
where $x_{pop}$ is the vector of coefficients, $b\in \Rm{m}{1}$ is the response vector, $A\in \Rm{m}{n}$ is a fixed design matrix and the noise vector $\xi\sim \mathcal{N}(0,\sigma^2 I)$. The full solution $x^*$ and the SAS solution $\hat x$ are seen as predictors of $x_{pop}$ which are evaluated by the mean square error (MSE)
\begin{equation}\label{eq:MSEdef}
   MSE(\hat x) = \avg{\|A(\hat x-x^*)\|^2_2}.
\end{equation}
Ma et al. \cite{ma15a} first analyzed the properties of SAS for uniform sampling and leverage score sampling in this model, looking both at the bias and the variance of the estimator. 

To compute a more accurate estimate of $x^*$ (in the optimization perspective) or $x_{pop}$ (in the statistical perspective), one can average $N$ SAS estimators, as suggested by \cite{sketched_ridgeReg,bartan,bartan2023distributed} and many others. 
Wang et al \cite{sketched_ridgeReg} perform a broad analysis for sketched ridge regression with averaging, 
both comparing the residual of the sketched solution with the optimal residual and analyzing the bias and variance in the statistical model \eqref{eq:model}. Defining $\bar x = \dfrac{1}{N}\sum_{i=1}^N \hat x_i$ the average of $N$ samples of SAS solutions, they prove that, for several sketching methods and sufficient sketch size,
\[
\|A\bar x - b\|_2^2 \leq \left(1+\dfrac{\eta}{N}  +\eta^2\right)\|Ax^*-b\|_2^2.
\]
with high probability.

Averaging also leads to new strategies to reduce the bias \cite{precond_viaDistrib, bartan2023distributed}. 

Furthermore, averaging several solutions was studied by Bartan and Pilanci \cite{bartan,bartan2023distributed} as a natural method for distributed computing, in the sense that several nodes compute the sketched solutions then the central node averages these solutions. In this way the distributed framework is also resilient to node failure. There are also SGD methods for asynchronous computing but their convergence speed depends on the conditioning of the problem \cite{bartan}.

In this paper we view the average of several SAS solutions as a Monte Carlo estimator, which inspired us to extend the analysis to a Multilevel Monte Carlo estimator.
Thus, we introduce an estimator called Multilevel Sketch-And-Solve (MLSAS) to approximate the solution $x^*$ and explore its advantages and drawbacks. 
The idea is inspired by the Multilevel Monte Carlo (MLMC) method introduced by Giles \cite{giles2008opre}, which was successfully applied as a variance reduction technique in stochastic path simulation, initially for financial pricing, and later in many more applications \cite{giles2015actanum}. Multilevel Monte Carlo uses different time stepping levels in the numerical solution of SDEs to correct the estimate obtained on the coarsest time stepping level, such that most samples of the MLMC estimator are computed very cheaply while also ensuring that the weak error (equivalent to our bias) of the MLMC estimator is that of the finest level.

In our context, instead of different levels of discretization of an SDE, we use different sketch sizes $s$ to  
estimate the solution $x^*$ of \eqref{eq:LSpb}. 

Multilevel ideas have already been applied in NLA to solve linear systems in \cite{acebron2020probabilistic} which uses Multilevel Monte Carlo to compute a Laplace transform. However, to the best of our knowledge this paper is the first that studies a Multilevel Monte Carlo approach for least squares problems.

\subsection{Main contribution: analysis of a Multilevel Sketch-And-Solve (MLSAS) estimator for the solution of the least squares problem}\label{sec:mainContrib}
We introduce the MLSAS estimator, which is a combination of sketch-and-solve and the Multilevel Monte Carlo paradigm for nested simulations \cite{giles2015actanum}. 
More precisely, let $L\in \mathbb{N}^*$ and for each level $\ell\in\{0, \ldots, L\}$ define the sketch size $s_{\ell}=2^{(\ell+1)} n$ (such that $n\leq s_\ell <m$) and a random sketch matrix $S_\ell\in \Rm{s_\ell}{m}$. 
We denote by $x^{(\ell)}$ the solution to problem \eqref{eq:SASpb} with $S=S_\ell$. To define the multilevel estimator, we also define the variables $\Delta x^{(\ell)} = x^{(\ell)}- x^{(\ell-1)} $ (with the convention $\Delta x^{(0)} = x^{(0)}$) where the coarser solution $x^{(\ell-1)}$ is based on a sketch matrix that is equal to $S_\ell(1:s_{\ell-1},:)$. 
We also consider the case where the matrix $S_\ell$ is split in two such that $S_\ell = \dfrac{1}{\sqrt{2}} \begin{bmatrix}
    S_a \\ S_b
\end{bmatrix}$ and $S_a,S_b$ are used to obtain two estimates $x^{(\ell-1)}_a,x^{(\ell-1)}_b$ that are averaged to get $x^{(\ell-1)}=\left(x^{(\ell-1)}_a+x^{(\ell-1)}_b\right)/2$. This is akin to \textit{antithetic variables} in the MLMC literature, which can lead to reduced variance in path simulation \cite{giles2014antithetic} and nested expectations \cite{giles2015actanum}. Therefore we also call this method antithetic variables in our work.

Then the MLSAS estimator of $x^*$ is written as
\begin{equation}\label{eq:mlmc_estim}
    x_{MLMC} = \sum_{\ell=0}^L N_\ell^{-1}\sum_{i=1}^{N_\ell} \Delta  x_i^{(\ell)},
\end{equation}
where $N_\ell$ is the number of random samples of $\Delta x^{(\ell)}$ used in the estimator. A crucial point is that the samples of $\Delta x^{(\ell)}_i$ are generated independently for all $\ell$ and $i$, unless we state otherwise (e.g. in Section \ref{sec:costAnalysis}).

To track the accuracy of the estimator, we use the Mean Squared Error (MSE) defined in \eqref{eq:MSEdef}. Denote by $V_\ell=\var{A\Delta x^{(\ell)}} = \avg{\|A\Delta x^{(\ell)}\|_2}- \|\avg{A\Delta x^{(\ell)}}\|_2$ the \textit{level variances}. The MSE can be decomposed into variance and bias contributions as 
\begin{equation}
       MSE(x_{MLMC}) =  \sum_{\ell=0}^L N_\ell^{-1} V_\ell + \|A\avg{x^{(L)}}-Ax^*\|^2_2.
\end{equation}
The variance per level $V_\ell$ is approximated as part of the estimation of $x^*$ and the bias can be estimated by an upper bound provided in \cite{bartan}. For uniform sampling, the bound on the bias is of the form
\[
\|A\avg{x^{(L)}}-Ax^*\|^2_2 \leq \sqrt{4 \eta \dfrac{n}{s_L}\|Ax^*-b\|_2^2 \max_i \|U_{i,:}\|_2^2}
\]
where $A =U\Sigma V^T$ is the thin SVD of $A$ and $U_{i,:}$ denotes the $i$-th row of the left singular vectors of $A$. Note that the above bound cannot be computed explicitly without solving the full least squares problem, however $\max_i \|U_{i,:}\|_2^2$ can be approximated using sketching \cite{drineas2012levSco}, and $\|Ax^*-b\|_2$ can be replaced by $(1+\eta)\|Ax^{(L)}-b\|_2$ for a computable bound.

The total computational cost of the MLSAS estimator is 
\begin{equation} 
    \mathcal{C}_{MLSAS} = \sum_{\ell=0}^L N_\ell C_\ell,
\end{equation}
where $C_\ell$ is the cost per sample of $\Delta x^{(\ell)}$.

The main motivation for applying an MLMC-like strategy to solve large-scale least squares problems is that for several common sketch types (uniform and leverage score sampling, SRHT, SRTT; see Section \ref{sec:BackgroundSketching} for their definition) the bias of the SAS estimate decreases with the sketch size $s$. 

In this paper, supported by numerical evidence and analysis, we show that 
\begin{enumerate}
    \item with antithetic variables, the variance $V_\ell$ decreases with a faster rate than without antithetic variables (see Section \ref{sec:errorAnalysis});
    \item 
    for a given desired total variance, the cost of the simple Monte Carlo estimator does not depend on the maximal level $L$, i.e. on the maximal sketch size (see Section \ref{sec:costAnalysis}); 
    \item the total cost of the MLSAS estimator is higher than the cost of simple Monte Carlo estimation (see Section \ref{sec:costAnalysis}).
\end{enumerate}
Thus a direct application of MLMC appears to be unattractive.

\subsection{Outline}
The remainder of the paper is organized as follows. 
First in Section \ref{sec:prelims} we summarize background on the sketching matrices that we are concerned with as well as the main ideas of Multilevel Monte Carlo as used in SDE applications. 
Then in Section \ref{sec:errorAnalysis} we analyze the variance of the estimator with and without antithetic variables. 
Finally, Section \ref{sec:costAnalysis} shows that the complexity of the multilevel SAS estimator cannot be lower than that of the simple Monte Carlo estimator.


\section{Preliminaries} \label{sec:prelims}
In this section we present the sketching methods that are relevant for the MLSAS estimation of the solution to \eqref{eq:LSpb}, then we summarize the core ideas of MLMC from \cite{giles2008opre} that are referred to in the discussion on the cost of the MLSAS method (Section \ref{sec:costAnalysis} of the paper).

\subsection{Subsampling and sketching matrices} \label{sec:BackgroundSketching}
Random sketching and subsampling allow one to reduce the size of the least squares problem, as discussed in the introduction. 
The key property to define sketching matrices is the ability to preserve the geometry of the original problem. This requirement is formalized through the subspace embedding property, which ensures that the action of the sketching matrix approximately preserves the Euclidean norm of all vectors in a given subspace.

Formally, let $ A \in \mathbb{R}^{m \times n} $ be a data matrix with rank $n$, and let $S \in \mathbb{R}^{s \times m}$ be a sketching matrix with $s \ll m$. We say that $S$ is a $(1 \pm \eta)$ subspace embedding for the column space of $A$ if, for all $x \in \mathbb{R}^n$,
\begin{equation} \label{eq:embedding}
(1 - \eta)\,\|Ax\|_2^2 \;\le\; \|SAx\|_2^2 \;\le\; (1 + \eta)\,\|Ax\|_2^2.
\end{equation}
This is a sufficient condition to obtain that
\begin{equation} \label{eq:embeddingMatrix}
    \bigl\| A^\top S^\top S A - A^\top A \bigr\|_2 \;\le\; \eta \,\|A\|_2^2,
\end{equation}
which highlights that the sketch approximately preserves the Gram matrix associated with $A$. The former can be proved by bounding $|x^T(A^\top S^\top S A - A^\top A)x|= |\|SAx\|_2^2- \|Ax\|_2^2|$ using the subspace embedding property. Equation \eqref{eq:embedding} ensures that the solution $\hat x$ to the sketched problem \eqref{eq:SASpb} satisfies 
\begin{equation}
    \|A\hat x-b\|_2^2 \leq \dfrac{1+\eta}{1-\eta} \|Ax^*-b\|_2^2,
\end{equation}
where $x^*$ is the solution of the full least squares problem.

Most theoretical accuracy guarantees in the sketching literature are established for Gaussian sketching matrices, owing to their strong concentration properties and analytical tractability. 
However, Gaussian matrices are dense and the cost of applying them is $\mathcal{O}(mns)$ which is avoided in practice by using more structured sketch matrices. Furthermore, as shown by Epperly in his blog \cite{epperly_blogBias}, for Gaussian sketching the SAS solution is unbiased. The same is not true of other, more efficient sketches.

We distinguish mainly several classes of sketch matrices:
(1) dense or sparse matrices and (2) data-aware or data oblivious sketches. Data-oblivious sketches are generated independently of the input matrix, whereas data-aware sketches depend on the data, for example for subsampling-based sketches,
in the definition of the probability of subsampling the rows. 
In this paper we focus on dense least squares problems and on sketching methods such as uniform sampling, leverage score sampling, and sketching matrices based on trigonometric and Hadamard transforms, which we present below. We also briefly introduce the Gaussian sketch, for completeness. 
Table \ref{tab:sketch_size} summarizes the sketch size that is necessary to ensure that the sketch matrix leads to the residual bound \eqref{eq:resBound}.
For uniform sketching, the sketch size depends on the coherence $\mu(A)$, which is defined as 
\[
\mu(A) \;=\; \frac{m}{n} \max_{1 \le i \le m} \|U_{i,:}\|_2^2,
\]
where $A = U \Sigma V^\top$ is a thin singular value decomposition and $U \in \mathbb{R}^{m \times n}$ has orthonormal columns.

\begin{table}[h!]
\centering
\caption{Sketch size sufficient to ensure a relative-error least-squares residual bound \eqref{eq:resBound} \cite[Theorem 1]{sketched_ridgeReg}. The cost of applying the sketch to a dense matrix $A\in\Rm{m}{n}$ is provided in the last column.}
\label{tab:sketch_size}
\begin{tabularx}{\textwidth}{l c X c}
\toprule
\textbf{Sketch type} & \textbf{Sketch size $s$} & \textbf{Assumptions} & Cost \\
\midrule
Gaussian
& $\mathcal{O}\!\left(\tfrac{n}{\eta}\right)$
& Oblivious sketch & $\mathcal{O}(smn)$ \\
SRTT / SRHT
& $\mathcal{O}\!\left(\tfrac{n \log n}{\eta}\right)$
& Oblivious sketch & $\mathcal{O}(mn\log n)$ \\ 

Leverage score sampling
& $\mathcal{O}\!\left(\tfrac{n \log n}{\eta}\right)$
& Sampling proportional to leverage scores & $\mathcal{O}(mn\log n)$ \\

Uniform subsampling
& $\mathcal{O}\!\left(\tfrac{\mu(A)\,n \log n}{\eta}\right)$
& Requires bounded coherence $\mu(A)$ & $\mathcal{O}(1)$ \\
\bottomrule
\end{tabularx}
\end{table}

\paragraph{Gaussian sketching}
Gaussian sketching corresponds to left multiplying $A$ by a matrix 
\[S=\dfrac{1}{\sqrt{s}}G\]
where all entries of $G$ are sampled independently from $\mathcal{N}(0,1)$. This sketching method is widely used for theoretical analysis, but is usually unattractive in practice because the cost of applying this sketch is $\mathcal{O}(mns)$, which is higher than for other sketching methods presented below.

\paragraph{Uniform sampling}
Uniform sampling corresponds to selecting $s$ rows of the data matrix $A \in \mathbb{R}^{m \times n}$ and the right-hand side $b \in \mathbb{R}^m$ uniformly at random. The associated sketching matrix $S \in \mathbb{R}^{s \times m}$ has exactly one nonzero entry per row, equal to $1$, indicating which row of $A$ (or $b$) is selected. When sampling is performed with replacement, the same row of $A$ may be selected multiple times.

To compensate for the reduction in sample size, the sketching matrix is rescaled by a factor $\sqrt{m/s}$, so that the expected second moment is preserved, i.e.,
\[
\mathbb{E}\!\left[ S^\top S \right] = I_m.
\]
This normalization ensures that, in expectation, the Euclidean norm of any vector is preserved under the action of the sketch, although uniform sampling does not in general satisfy the subspace embedding property unless the rows of $A$ have approximately uniform leverage scores. The residual bound \eqref{eq:resBound} holds w.h.p. for this sketching method for $s = \mathcal{O}(\mu(A) n \log n / \eta)$ but can fail if the coherence $\mu(A)$ of the matrix $A$ is high \cite[Sec. 6.1.1]{murray2023randomized}.

\paragraph{Leverage score sampling}
Leverage score sampling is a data-aware sketching strategy that selects rows of the matrix $A \in \mathbb{R}^{m \times n}$ according to their statistical importance \cite[Sec. 6]{murray2023randomized}. Let $A = U \Sigma V^\top$ be a thin singular value decomposition of $A$, where $U \in \mathbb{R}^{m \times n}$ has orthonormal columns. The leverage score associated with the $i$-th row of $A$ is defined as
\begin{equation}\label{eq:LevScoresDef}
    \ell_i = \|U_{i,:}\|_2^2, \qquad i = 1, \dots, s,
\end{equation}
and satisfies $ \sum_{i=1}^m \ell_i = n $.

In leverage score sampling, rows are selected independently with probabilities proportional to their leverage scores. The sketching matrix $S \in \mathbb{R}^{s \times m}$ has one nonzero entry per row, corresponding to a selected row index $i$, and is rescaled by $1/\sqrt{s\,p_i}$, where $p_i = \ell_i / n$ denotes the sampling probability. This normalization ensures that $\mathbb{E}[S^\top S] = I_m $. Therefore the sketching matrix $S$ can be written as 
\[S = \dfrac{1}{\sqrt{s}} DR \]
where $R$ represents the selection of the rows and $D$ contains the scaling by the leverage scores.

Leverage score sampling requires that $s = \mathcal{O}(n \log n / \eta)$ to ensure that the residual guarantee \eqref{eq:resBound} holds with high probability \cite{sketched_ridgeReg}. 

In practice, since computing the exact leverage scores would be almost as expensive as solving the original least squares problems, one uses sketching to cheaply approximate the leverage scores, as studied in \cite{drineas2012levSco}.

\paragraph{SRHT and SRTT sketching}
The Subsampled Random Hadamard Transform (SRHT) and the Subsampled Random Trigonometric Transform (SRTT) are structured random sketching methods that combine randomization with fast orthogonal transforms. Both sketches can be written in the form
\[
S = \sqrt{\frac{m}{s}}\, R T D,
\]
where \(D \in \mathbb{R}^{m \times m}\) is a random diagonal matrix with independent Rademacher entries ($1$ or $-1$ with probability $0.5$), \(T \in \mathbb{R}^{m \times m}\) is an orthogonal transform, and \(R \in \mathbb{R}^{s \times m}\) is a uniform sampling matrix that selects \(s\) rows without additional rescaling. In the SRHT case, \(T\) is a (normalized) Hadamard matrix, while for the SRTT, \(T\) corresponds to an orthogonal trigonometric transform. 

The random sign matrix $D$ serves to randomize the input, while the orthogonal transform $T$ spreads the energy of the rows of $A$ approximately uniformly, therefore applying the transformation $TD$ to a matrix reduces its coherence. As a result, the subsequent uniform sampling step $R$ is effective, and the combined sketch ensures that the error bound \eqref{eq:resBound} is met with high probability for $s=\mathcal{O}(n\log n/\eta)$ \cite{tropp2011improved}
regardless of the properties of the matrix $A$. 

By construction of this sketching technique we only need to apply the matrix $TD$ once then generate different samples by varying the rows selected by $R$, therefore the cost of applying $TD$ is not included in the cost per sample.

\subsection{Multilevel Monte Carlo for SDEs}\label{sec:BackgroundMLMC}
Multilevel Monte Carlo (MLMC) was introduced by Giles \cite{giles2008opre} as a variance reduction method for path simulation. The most common application of MLMC is the construction of numerical methods for financial options pricing, and numerous variants and extensions of the method are presented in \cite{giles2015actanum}.

The principle of MLMC  \cite{giles2008opre} is as follows: 
let $X_t$ be a variable that follows a stochastic differential equation of the form
\[dX_t = a(t,X_t)dt + b(t, X_t) dW_t.\]
Suppose we want to estimate the expectation $\mathbb{E}[f(X_T)]$ of a function that depends on the solution $X_t$ at time $T$. Let $P=f(X_T)$. The classic Monte Carlo method then consists of discretizing the time interval $[0,T]$ into $n$ steps of length $h$ (so $t=kh$ with $0\leq k\leq n$) and simulating $N$ paths $\widehat X_k(\omega)$ (where $\omega=1, \ldots, N$ are the random events that occur) that approximate the path of $X_t$, estimating the expectation by the mean of the samples of $\widehat P(\omega) =f(\widehat X_n(\omega))$:
\[
\avg{P} \approx \dfrac{1}{N}\sum_{\omega=1}^N \widehat P(\omega) =:\widehat P_{MC} .
\]
Assuming that the samples are independent, the variance of the above estimator is $\mathbb{V}[\widehat P]/N$. Therefore, if we want to have a standard deviation $\varepsilon$ and a weak error $\varepsilon$, we must simulate $N = \mathcal{O}(\varepsilon^{-2})$ paths with cost $\mathcal{O}(h^{-1})=\mathcal{O}(\varepsilon^{-1})$, so the total cost of the estimation is $\mathcal{O}(\varepsilon^{-3})$ \cite{giles2008opre}.

To reduce this cost, in MLMC one calculates paths with different time steps $h_\ell$ (for example, $h_\ell = 2^{-\ell}$, for $\ell=0, \ldots , L$) to approximate $P^{(\ell)}$, then forms variables $\Delta P^{(\ell)} = P^{(\ell)} - P^{(\ell-1)}$. The terms in $\Delta P^{(\ell)}$ depend on the same Gaussian increments that approximate Brownian motion $dW$ so that $\Delta P^{(\ell)}$ can be seen as a term that corrects the bias of $P^{(\ell-1)}$. By combining these terms for levels $0$ to $L$ and defining $\Delta P^{(0)} = P^{(0)}$, we obtain the following identity:
\[\avg{P^{(L)}} = \sum_{\ell=0}^L \avg{\Delta P^{(\ell)}}.\]

We then define the MLMC estimator of $\avg{P}$ by estimating each term in the above sum by the average of $N_\ell$ samples, which gives:
\[\widehat P_{MLMC} =\sum_{\ell=0}^L N_\ell^{-1} \sum_{i=1}^{N_\ell} \Delta P^{(\ell)}(\omega_i).\]

Assuming that all samples of $\Delta P^{(\ell)}$ for all $\ell$ are calculated independently, the total variance of the MLMC estimator is therefore
\[
V_{tot} = \sum_{\ell=0}^{L} \dfrac{V_\ell}{N_\ell}
\]
where $V_\ell = \var{\Delta P^{(\ell)}}$ is the variance for each level. Defining $C_\ell$ to be the cost of calculating a sample of $\Delta P^{(\ell)}$, the total cost of the estimator is
\[
C_{tot} = \sum_{\ell=0}^{L} N_\ell C_\ell.
\]
The advantage of MLMC over the classic Monte Carlo estimator is that it provides an estimate whose bias is equal to that of the finest level and that the variances $V_\ell$ are much smaller than the terms $\var{P^{(\ell)}}$, so for the same value $\varepsilon$ of the desired standard deviation, the MLMC estimator requires fewer expensive calculations.

Indeed, by using a Lagrange multiplier $\lambda^2$ and considering the number of samples $N_\ell$ for each level as a decision variable, we can minimize the cost $C_{tot}$ under the constraint $V_{tot}=\varepsilon^2$ \cite{giles2008opre}. At the optimum, we have $N_\ell = \left \lceil \varepsilon^{-2}\sqrt{V_\ell/C_\ell} \sum_{k=0}^L \sqrt{V_k C_k} \right \rceil$ and the total cost of the estimator is
\begin{equation}\label{eq:costMLMC}
    C_{MLMC} = \varepsilon^{-2} \left(\sum_{\ell=0}^L \sqrt{V_\ell C_\ell}\right)^2.
\end{equation}
Whereas for the same bias with standard Monte Carlo estimation, the number of samples required is $N=\left \lceil V_L^{MC}/\varepsilon^2\right \rceil$ and the total cost is
\begin{equation}\label{eq:costMC}
    C_{MC} = \varepsilon^{-2}V_L^{MC} C_L^{MC}.
\end{equation}

This leads to the MLMC algorithm \ref{alg:MLMCalgo} for estimation of $\avg{P}$. The number of levels $L$ is defined such that the weak error $\avg{P^{(L)}-P}$, which is considered as the bias of the estimator, is smaller than a specified tolerance (typically $\varepsilon/2$). In practice, the number of levels is updated and new samples are computed iteratively until the estimated total MSE reaches the desired value. Assuming that $\avg{P_\ell - P_{\ell -1}} \propto 2^{-\alpha \ell}$, the test for weak convergence is $|\avg{P_L-P_{L-1}}|/(2^\alpha-1)\leq tol$, and the variances $V_\ell$ are estimated using the samples computed in the main loop. 
Heuristically, the algorithm converges and reaches the desired accuracy as long as the estimators $\Delta P^{(\ell)}$ satisfy the assumptions of \cite[Theorem 3.1]{giles2008opre}. These assumptions hold in many practical cases, for example when using the Euler-Maruyama numerical scheme and a Lipschitz payoff function. The computational savings in MLMC come from the fact that most samples are computed on the coarsest level, therefore, on top of a lower total computational cost, the MLMC framework offers high potential for parallelism.

\begin{algorithm}[h]
\caption{Multilevel Monte Carlo algorithm}\label{alg:MLMCalgo}
\begin{algorithmic}[1]
\STATE Input: parameters describing the SDE, initial maximal level $L$, number of samples to start with $N_0$ (e.g. $N_0=10^4$) and desired RMSE $\varepsilon$.
\STATE Output: estimate of $\avg{P}$
\STATE generate $N_0$ samples on levels $\ell=0,\ldots,L$ and estimate the variances $V_\ell$
\WHILE{$\sum_\ell N_\ell^{-1}V_\ell>\varepsilon^2/2$}
\STATE calculate optimal number of samples $N_\ell$
\STATE generate extra samples at each level $\ell$ as required depending on $N_\ell$
\STATE test for weak convergence: if $|\avg{P_L-P_{L-1}}|/(2^\alpha-1)> \varepsilon^2/2$, add a new level $L:=L+1$
\ENDWHILE
\end{algorithmic}
\end{algorithm}


\section{Error analysis of the multilevel SAS estimator} \label{sec:errorAnalysis}
\begin{figure}[htb]
    \centering
    \includegraphics[width=0.8\linewidth]{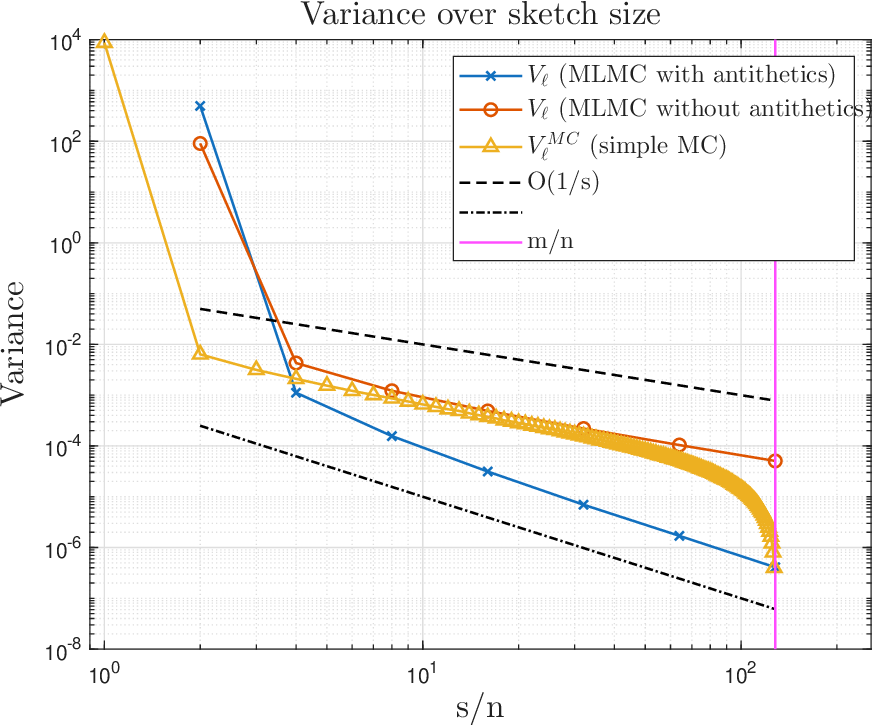}
    \caption{Comparison of the variances $V_\ell=\var{Ax^{(\ell)}-Ax^{(\ell-1)}}$, $V_\ell=\var{Ax^{(\ell)}-A(x_a^{(\ell-1)}+x_b^{(\ell-1)})/2}$ and $V^{MC}_\ell$ for different levels. Here $A$ is incoherent ($U\in \Rm{m}{n}$, the matrix of left singular vectors of $A$ is Haar distributed) and we used uniform sampling with replacement. Other sampling or sketching methods lead to similar results.}
    \label{fig:var_antit}
\end{figure}

In this section we analyze the variance of MLSAS (with and without antithetics). 
We begin with an experiment to illustrate the variance reduction offered by the multilevel method. 
We take a matrix $A$ with dimensions $m=6400$ by $n=50$ using the MATLAB command $\texttt{gallery(‘randsvd’,[m, n],1e2)}$, which generates a banded matrix with geometrically distributed singular values, with condition number $\kappa(A)=100$ and Haar distributed singular vectors, and we define the vector $b$ by $\texttt{b=A*randn(n,1) + randn(m,1)*1e-3}$. 
We generate 1000 samples for each level $\ell$, and calculate the variance of $x^{(\ell)}$ and $\Delta x^{(\ell)}$, using uniform sampling with replacement for different sketch sizes. We display the variance over sketch size in Figure \ref{fig:var_antit}.

We observe that the use of antithetic variables changes the slope of the variance compared to simple Monte Carlo and Multilevel Monte Carlo without antithetics. It is in fact necessary to use antithetic variables if we want to reduce the variance compared to the Monte Carlo estimator. We thus observe that $V_\ell \propto 4^{-\ell} \propto 1/s_\ell^2 $. The difference in trend observed with or without antithetic variables can be explained by a quick analysis of the term $A\Delta x^{(\ell)} $, which we detail below.

We obtained similar results using SRTT sketching and leverage score sampling, which are applicable even when the matrix $A$ has a higher coherence. With SRTT, the application of $TD$ is known to reduce the coherence, so one can compute $TDA$ at the outset, then turn to subsampling-based sketches to the problem $\min_x\|(TD)(Ax-b)\|_2$. 

For leverage score sampling, if the right hand side $b$ has large components in the rows with small leverage scores, we advocate to use leverage scores computed from the SVD of $[A,\: b]$ instead of the classical definition \eqref{eq:LevScoresDef}. In addition, we also used weights $\mbox{diag}(1/\sqrt{\ell_i})$ in solving the sketched least squares problem with leverage score sampling, as we observed it made the estimated variances more robust and reduced them by a considerable constant factor. As proved in \cite{jones2025subapsnap}, this version of the leverage scores notably ensures that the quantity $\|A\hat x-b\|_2/\|Ax^*-b\|_2$ is of order $\mathcal{O}(1)$, instead of $\mathcal{O}(1/\sigma_{min}(SQ))$ (where $Q$ is from the QR decomposition of $A$). In our experiments this definition of the leverage scores led to the same variance trends as the ones observed in Figure \ref{fig:var_antit}. On the other hand, using leverage scores based only on $A$ increased the variances by a large factor $10^4$ and reduced the slope of the variance obtained with MLMC with antithetics.

\subsection{Variance of the MLSAS estimator}
We now explain the trend of the level variance $V_\ell=\var{ A \Delta x^{(\ell)}}$ analytically and experimentally. The variance can be decomposed as 
\[
V_\ell = \avg{\|A\Delta x^{(\ell)}\|_2^2} - \|A\avg{\Delta x^{(\ell)}}\|_2^2
\]
and we define $\alpha_\ell = \avg{\|A\Delta x^{(\ell)}\|_2^2}$ and $\beta_\ell = \|A\avg{x^{(\ell)}}\|_2^2$. In our experiments we observed that $\beta_\ell \ll \alpha_\ell$, therefore in the analysis we focus on the term $\alpha_\ell$, and assume $V_\ell \approx \alpha_\ell$.

\paragraph{Expression of $\alpha_\ell$ with antithetic variables}
Suppose that the sampling matrices $S_\ell, S_a,S_b$ are fixed. Let $A=U\Sigma V^T$ be the SVD; we get by definition of the SAS solution that (assuming $S_\ell A$ is full rank)
\begin{align*}
x^{(\ell)}  &= (S_\ell A)^\dagger S_\ell b \\
&=(A^TS_\ell^TS_\ell A)^{-1}A^T S_\ell^T S_\ell b \\
&= V\Sigma^{-1} (U^T S_\ell^T S_\ell U)^{-1} \Sigma ^{-1}V^T V\Sigma U^T S_\ell^T S_\ell b \\ 
&= V\Sigma^{-1} (U^T S_\ell^T S_\ell U)^{-1} U^T S_\ell^T S_\ell b
\end{align*}
To simplify the notation, for the indices $i=a,b,\ell$ we define $H_i = U^TS_i^TS_iU $ and $ P_i = S_i^TS_i$. Multiplying $x^{(\ell)}$ by $A$ on the left, we obtain 
\[
A x^{(\ell)}  = U(U^T S_\ell^T S_\ell U)^{-1} U^T S_\ell^T S_\ell b = U H_\ell^{-1}U^TP_\ell b.
\]
We notice that
\[
S_\ell^TS_\ell = \dfrac{1}{2} \begin{bmatrix} S_a^T & S_b^T\end{bmatrix} \begin{bmatrix} S_a \\S_b \end{bmatrix} = \dfrac{1}{2}(S_a^TS_a + S_b^TS_b).
\]
and therefore $H_\ell = \dfrac{1}{2} (H_a+H_b)$ as well as $P_\ell = \dfrac{1}{2}(P_a+P_b)$.
Therefore, using antithetic variables leads to
\begin{align*}
    A\Delta x^{(\ell)}  &= U \left((H_a+H_b)^{-1} U^T (P_a+P_b) -\dfrac{1}{2}(H_a^{-1}U^TP_a + H_b^{-1} U^T P_b) \right)b \\
    &= U (H_a+H_b)^{-1} \left(U^T (P_a+P_b) - \dfrac{1}{2}U^T (P_a+P_b)- \dfrac{1}{2}(H_b H_a^{-1}U^TP_a + H_a H_b^{-1} U^T P_b) \right)b \\
    &=\dfrac{1}{2} U (H_a+H_b)^{-1} \left(U^T (P_a+P_b)- (H_b H_a^{-1}U^TP_a + H_a H_b^{-1} U^T P_b) \right)b \\
    &= \dfrac{1}{2} U (H_a+H_b)^{-1} \left((I-H_b H_a^{-1})U^TP_a + (I-H_a H_b^{-1}) U^T P_b \right)b \\
    &= \dfrac{1}{2} U (H_a+H_b)^{-1} \left((H_a-H_b )H_a^{-1}U^TP_a + (H_b-H_a) H_b^{-1} U^T P_b) \right)b \\
    &= \dfrac{1}{2} U (H_a+H_b)^{-1} (H_a-H_b )(H_a^{-1}U^TP_a - H_b^{-1}U^T P_b)b \\
    &= \dfrac{1}{2} U (H_a+H_b)^{-1} (H_a-H_b )\Sigma V^T (x_a^{(\ell-1)} - x_b^{(\ell-1)})
\end{align*}

The intuition that we will check later on is therefore that the factors $''a-b''$ have smaller norm than their analogue without antithetic variables.

\paragraph{Expression of $\alpha_\ell$ without antithetic variables}
    On the other hand, without antithetic variables one can simply consider that $\Delta x^{(\ell)} = x^{(\ell)} - x_a^{(\ell-1)}$, with $x_a^{(\ell-1)}$ defined as above. Then, with the same notation as above we have
\begin{align*}
    A\Delta x^{(\ell)} &= U \left((H_a+H_b)^{-1} U^T (P_a+P_b) -H_a^{-1}U^TP_a \right)b \\
    &= U (H_a+H_b)^{-1} \left(U^T (P_a+P_b) -(H_a+H_b)H_a^{-1}U^TP_a \right)b  \\
    &= U (H_a+H_b)^{-1} \left(U^T P_b -H_bH_a^{-1}U^TP_a \right)b  \\
    &= U (H_a+H_b)^{-1}H_b \left(H_b^{-1}U^T P_b -H_a^{-1}U^TP_a\right)b \\
    &= U (H_a+H_b)^{-1}H_b \Sigma V^T (x_b^{(\ell-1)} -x_a^{(\ell-1)})
\end{align*}
Here instead of the factor $H_a-H_b$ we obtain $H_b$, which explains the difference in the trend of $V_\ell$ with or without antithetic variables. Hence we have proved the proposition below. Note that the two decompositions of $V_\ell$ above hold for each of the sketching or subsampling methods from the Section \ref{sec:BackgroundSketching}, since we only used the general definition of $S_\ell$ as the concatenation of $S_a,S_b$. 
\begin{prop}
    Let $x_a^{(\ell-1)}, x_b^{(\ell-1)}$ obtained from solving the sketched problem \eqref{eq:SASpb} with the sketch matrices $S=S_a,S_b$ respectively, and $x^{(\ell)}$ the corresponding finer solution obtained from solving \eqref{eq:SASpb} with $S=\dfrac{1}{\sqrt{2}}\begin{bmatrix}
        S_a \\ S_b
    \end{bmatrix}$.
    For indices $i=a,b,\ell$ note $H_i = U^TS_i^TS_iU $ and $ P_i = S_i^TS_i$. 
    
    If we define $\Delta x^{(\ell)}= x^{(\ell)}-\dfrac{1}{2}(x_a^{(\ell-1)}+x_b^{(\ell-1)})$ (ie. we use antithetic variables), we have
    \begin{equation}
        A\Delta x^{(\ell)} = \dfrac{1}{2} U (H_a+H_b)^{-1} (H_a-H_b )\Sigma V^T (x_a^{(\ell-1)} - x_b^{(\ell-1)}).
    \end{equation}
    On the other hand, if we define $\Delta x^{(\ell)}= x^{(\ell)}-x_a^{(\ell-1)}$ (ie. no use of antithetic variables) then
    \begin{equation}
        A\Delta x^{(\ell)} = U (H_a+H_b)^{-1}H_b \Sigma V^T (x_b^{(\ell-1)} -x_a^{(\ell-1)}).
    \end{equation}
\end{prop}

\begin{figure}[htb]
    \centering
    \includegraphics[width=0.7\linewidth]{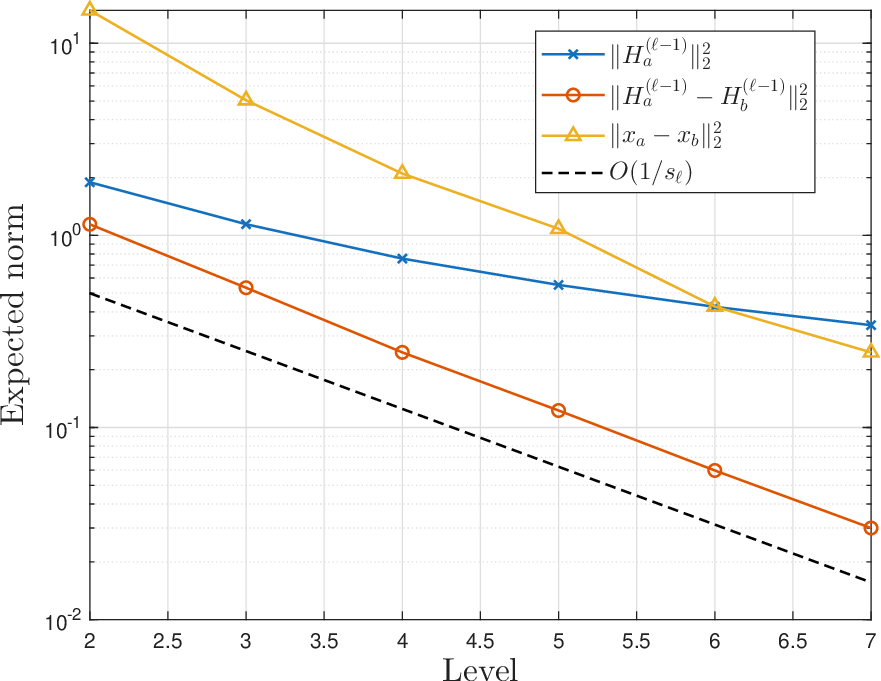}
    \caption{Average over 100 runs of spectral norm of the factors of  $\|A\Delta x^{(\ell)}\|^2$, for sketch size $s_\ell = n2^\ell$ on level $\ell$.}
    \label{fig:Hnorm}
\end{figure}

To verify that the decomposition above explains the difference in trends with or without antithetics, let us consider a numerical example.
We used the same matrix $A$ and vector $b$ as before and calculated the squared norms of each factor appearing in $A\Delta x^{(\ell)}$ for $N=100$ samples, then took the average of these norms. The results are shown in Figure \ref{fig:Hnorm}.

Table \ref{tab:slopesHnorm} compares the trends observed in Figure \ref{fig:Hnorm} for each factor of $A\Delta x^{(\ell)}$. 
By taking the product of the observed slopes (table \ref{tab:slopesHnorm}), we find the curves of \Cref{fig:var_antit}, which shows that this decomposition of $\|A\Delta x^{(\ell)}\|^2$ explains well the slopes observed in $V_\ell$ (\Cref{fig:var_antit}). In fact, using subspace embedding probabilistic guarantees \cite[Prop. 6.1.1]{murray2023randomized}, we can prove that $\|H_a-H_b\|^2 = \mathcal{O}(1/s)$ for uniform and leverage score sampling. For a distortion parameter $\eta \in [0,1]$, if
\[
2n\left( \dfrac{e^\eta}{(1+\eta)^{1+\eta}}\right)^{s/n} \geq \delta
\]
then the subspace embedding property \eqref{eq:embedding} fails with probability at most $\delta$ \cite[Prop. 6.1.1]{murray2023randomized}. Taking the logarithm gives 
\begin{equation} \label{eq:bdH1H2}
\dfrac{s}{n} \left[ \epsilon - (1+\epsilon)\ln(1+\epsilon) \right] \ge \ln(\delta/2n).
\end{equation}
In the regime where $s\gg n$, the distortion parameter $\eta$ is small, therefore we can use the Taylor expansion $\ln(1+\eta) \approx \eta - \frac{\eta^2}{2}$. Plugging this in \eqref{eq:bdH1H2} leads to 
$\dfrac{s}{n} \left( -\frac{\eta^2}{2} \right) \geq \ln(\delta/2n)$, and after rearranging to bound $\eta$ we obtain
\[
\eta \leq \sqrt{\dfrac{-2n\ln(\delta/2n)}{s}} = \mathcal{O}(1/\sqrt{s}).
\]
Finally, the matrix $S$ is a subspace embedding for $\mbox{range}(A)$ so it is a subspace embedding for $U$ (from the SVD of $A$), therefore
\begin{align*}
    \|H_a-H_b\|_2^2 &\leq (\|H_a-I\|_2 +\|H_b-I\|_2)^2 \\
    & \leq 2(\|H_a-I\|_2^2 +\|H_b-I\|_2^2) \\
    &= 2(\|(S_aU)^TS_aU -I\|_2^2 +\|(S_bU)^TS_bU -I\|_2^2) \\
    &\leq 4\eta^2\\  
    &= \mathcal{O}(1/s),
\end{align*} 
where the third line follows from \eqref{eq:embeddingMatrix}. This proves that $\|H_a-H_b\|_2^2=\mathcal{O}(1/s)$ with high probability.

\begin{table}[H]
    \centering
\caption{Asymptotic scaling of the spectral norm of the factors appearing in the expression of $A\Delta x^\ell$, according to what we observe from \Cref{fig:Hnorm}.}
    \label{tab:slopesHnorm}    
    \begin{tabular}{|p{4cm}|p{4cm}|} 
    \hline
         With antithetics & Without antithetics  \\
         \hline
        \multicolumn{2}{|c|}{$\|H_\ell^{-1}\|^2 = \mathcal{O}(1)$} \\ 
        \multicolumn{2}{|c|}{$\|x_a^{(\ell-1)}-x_b^{(\ell-1)}\|^2 = \mathcal{O}(1/s)$} \\
        \hline
        $\|H_a-H_b\|^2 = \mathcal{O}(1/s)$ & $\|H_b\|^2 = \mathcal{O}(1)$ \\
        $\|A\Delta x^\ell\|^2 = \mathcal{O}(1/s^2)$ & $\|A\Delta x^\ell\|^2 = \mathcal{O}(1/s)$ \\
        \hline
    \end{tabular}    
\end{table}

\subsection{Variance of the simple Monte Carlo estimator}\label{sec:varSimpleMC}
For the simple Monte Carlo estimator, the observed variance in Figure \ref{fig:var_antit} is usually $\mathcal{O}(1/s)$. The theoretical result in \cite[Lemma 1]{bartan} predicts this slope for Gaussian sketching with $s>n+1$,
\[
V_{MC}=\avg{\|A\hat x -Ax^*\|_2^2}\leq \dfrac{n}{s-n-1} \|Ax^*-b\|_2^2.
\]


\section{Cost of the multilevel estimator} \label{sec:costAnalysis}
We now analyze the cost of the MLSAS estimator and compare it to that of simple Monte Carlo. We show that, for all sketching methods from Section \ref{sec:BackgroundSketching} the overall cost of the simple Monte Carlo approach for a desired total variance of $\varepsilon^2$ is constant over level. Then deriving the total cost of the MLSAS estimation from the MLMC theory (see Section \ref{sec:BackgroundMLMC}) we show that the MLSAS framework cannot improve on the simple Monte Carlo framework in this application, including if we use a method to reduce the cost per sample $C_\ell$ on levels $\ell= 2, \ldots, L$.

\paragraph{Classical definition of the cost}
To evaluate the cost of the MLMC estimator, we would naively define the computational cost $C_\ell$ of a sample of $\Delta x^\ell$ (with antithetic variables) as $cost(x^{(\ell)}) + 2\times cost(x^{(\ell-1)})$. Using the Householder QR algorithm \cite[Sect. 5.2.2]{golub2013matrix}, each sample $x^{(\ell)}$ costs $2s_\ell n^2-\dfrac{2}{3}n^3+\mathcal{O}(mn)$ flops, therefore we would get $C_\ell = 2s_\ell n^2 + 4s_{\ell-1} n^2-2n^3=4s_\ell n^2-2n^3$. 
However, instead of this, we can use the calculations already performed on the antithetic variables $x_a^{(\ell-1)},x_b^{(\ell-1)}$ to obtain the finer variable $x^{(\ell)}$. Indeed, knowing the QR decompositions calculated for $x_a^{(\ell-1)},x_b^{(\ell-1)}$, i.e., that $S_a A= Q_aR_a$ and $S_b A= Q_bR_b$, we can note that 
\[S_\ell A = \dfrac{1}{\sqrt{2}}\begin{bmatrix}
    Q_a & 0 \\ 0 & Q_b
\end{bmatrix} \begin{bmatrix}
    R_a \\ R_b
\end{bmatrix}
\]
so it suffices to compute the QR factorization of $\begin{bmatrix}
    R_a \\ R_b
\end{bmatrix}$ to obtain $x^\ell$, which costs $\dfrac{10}{3}n^3$. 
Therefore the cost per sample of $\Delta x^{(\ell)}$, for $\ell\geq 1$, is $C_\ell = 4s_{\ell-1}n^2-\dfrac{4}{3}n^3 +\dfrac{10}{3}n^3=2s_{\ell}n^2+2n^3$. 
However our experiments show that with this expression of the cost the total cost \eqref{eq:costMLMC} is higher than \eqref{eq:costMC}.

\begin{figure}[htb]
\centering
\begin{tikzpicture}[>=Stealth, every node/.style={font=\Large$#1$}]

\node (Nl)   at (-4,0) { \text{Level } $\ell$ };
\node (Nl1)  at (-4,-2) { \text{Level } $\ell-1$ };
\node (Nl2)  at (-4,-4) { \text{Level } $\ell-2$ };

\node (Xl)   at (1.5,0)   { $x^{(\ell)}$ };

\node (Xa)   at (0,-2)   { $x_a^{(\ell-1)}$ };
\node (Xb)   at (3,-2)   { $x_b^{(\ell-1)}$ };

\node (Xaa)  at (-1,-4)  { $x_{aa}^{(\ell-2)}$ };
\node (Xab)  at (1,-4)  { $x_{ab}^{(\ell-2)}$ };
\node (Xba)  at (3,-4)  { $x_{ba}^{(\ell-2)}$ };
\node (Xbb)  at (5,-4)  { $x_{bb}^{(\ell-2)}$ };

\draw[<-] (Xl) -- (Xa);
\draw[<-] (Xl) -- (Xb);
\draw[<-] (Xa) -- (Xaa);
\draw[<-] (Xa) -- (Xaa);
\draw[<-] (Xa) -- (Xab);
\draw[<-] (Xb) -- (Xba);
\draw[<-] (Xb) -- (Xbb);

\node (dxl) at (7,0)   { $\Delta x^{(\ell)} $};
\node (dxl1) at (7,-2)  { $\left(\Delta x^{(\ell-1)}\right)$ };
\node (dxl2) at (7,-4)  { $\Delta x^{(\ell-2)}$ };

\draw[->] (Xl) -- (dxl);
\draw[->] (Xb) -- (dxl);
\draw[->] (Xbb) -- (dxl2);

\end{tikzpicture}
\caption{Method of reusing samples from previous levels to compute the samples on next levels: four samples ($x_{aa},x_{ab}, x_{ba},x_{bb}$) from level $\ell-2$ are used to compute  $\Delta x^\ell$. They are also used as the finer sample in $\Delta x^{\ell-2}$ but not used in $\Delta x^{\ell-1}$ so that they don't cancel out when the samples are summed in \eqref{eq:mlmc_estim}. In this figure arrows indicate in which quantities each SAS sample is used. }

\label{fig:skip_level}
\end{figure}
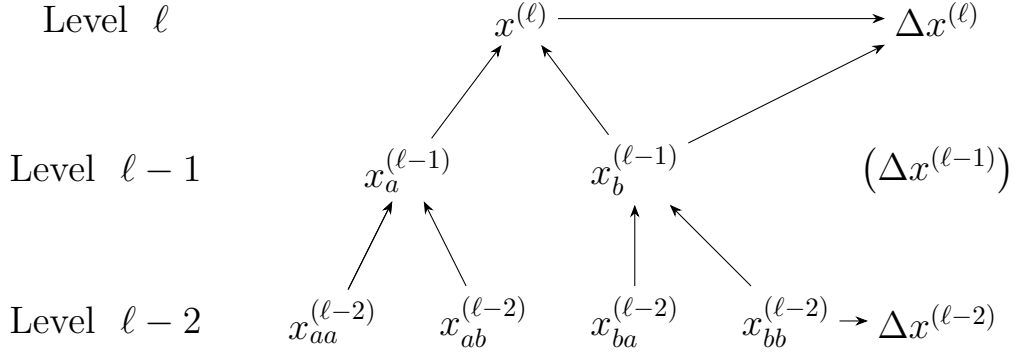

\paragraph{Attempt to reduce the cost per level} 
To reduce costs, we can think of reusing samples from lower levels to calculate those from higher levels, i.e., samples from the first levels are “recycled” in subsequent levels. To avoid that the samples from two consecutive levels cancel out in \eqref{eq:mlmc_estim}, we suggest to use the samples from level 0 for even levels and those from level 1 for odd levels, so that each sample appears either exclusively as an $x^{(\ell)}$ or in an antithetic variable (exclusive “or”), while being based on samples already calculated at lower levels; see diagram \Cref{fig:skip_level} for an illustration.

Using the formula of the optimal number of samples $N_\ell$, we can show that for all $\ell=2k+1$, $N_\ell \leq 2^{-\ell} N_1$ and for all $\ell=2k$, $N_\ell \leq 2^{-\ell} N_0$, which means that no extra samples on levels $\ell=0,1$ need to be generated in order to compute the samples on higher levels.  

However an important drawback of reusing samples from previous levels is that samples from different levels are \textit{not} independent, contrary to the assumption generally made in MLMC applications. 
To remain rigorous, covariance terms should theoretically be included in the expression of total variance. Since it would be simpler to neglect covariances, we wanted to see if they were significant in practice. 
In an experiment, we calculated $N_0$ (resp. $N_1$) samples at levels $\ell=0,1$ and calculated $N_0/4$ corresponding samples at level $\ell=2$ (resp. $N_1/4$ samples at level $\ell=3$). We then calculated their covariances (using the formula $Cov(X,Y) =\avg{XY^T}- \avg{X}\avg{Y}^T$) and found that: 
\begin{itemize}
    \item looking at the correlations of $\Delta x^l$:  $\rho(\ell=0,\ell=2)= 1.7e-2$  and  $\rho(\ell=1,\ell=3)= 5e-5$.
    \item looking at the correlations of the $A\Delta x^l$: $\rho(\ell=0,\ell=2)= 4.5e-7$  and  $\rho(\ell=1,\ell=3)= 6e-5$.
\end{itemize}
Therefore we can consider that the levels are weakly correlated. 
With this change, and taking $s_0=4n$, the level costs per sample are $C_0=\left(8-\dfrac{2}{3}\right)n^3\approx 7n^3$, $C_1=18n^3 $ and, for $\ell\geq 2$, $C_\ell = \dfrac{10}{3}n^2 \approx 3n^2$. 
Despite this, the cost of the MLSAS estimation is still higher than that of simple Monte Carlo. Indeed assuming that the variance $V_\ell=4^{-\ell}V_0$, where $V_0=\var{\Delta x^{(0)}}=\var{x^{(0)}}$, the cost of simple Monte Carlo, after choosing the right level $L$ depending on the tolerance one has on the bias (see MSE splitting \eqref{eq:MSEdef}), is 
\begin{align*}
    \varepsilon^2 \times \mathcal{C}_{MC} &= V_L^{MC} C_L^{MC} \\
    &\approx 2^{-L}V_0 \times 2^{L+3} n^3 = 8V_0 n^3
\end{align*}
On the other hand, from \eqref{eq:costMLMC} we have
\begin{align*}
    \varepsilon^2 \times \mathcal{C}_{MLSAS}  &= \left( \sum_{\ell=0}^L  \sqrt{V_\ell C_\ell}\right)^2 \\
    &\approx \left( \sqrt{7n^3V_0}+\sqrt{18n^3V_0/4}+ \sqrt{3V_0 n^3}\sum_{\ell=2}^L  2^{-\ell}\right)^2 \\
    &= n^3 V_0 \left( \sqrt{7} +\sqrt{18/4}+\dfrac{\sqrt{3}}{2}(1-2^{-L+1}) \right)^2 \\
    &> 8n^3 V_0
\end{align*}
Hence, despite the suggested reduction in cost per sample at higher levels, the multilevel estimation is still more expensive than the simple Monte Carlo estimation. We summarize this result in the proposition below.
\begin{prop}
    For a fixed tolerance on the total MSE, the Multilevel Sketch-And-Solve estimator \eqref{eq:mlmc_estim} is more expensive to evaluate than the simple Monte Carlo estimator $\dfrac{1}{N}\sum_{i=1}^N \hat x_i $, where each $\hat x_i $ is a sample computed with SAS.
\end{prop}

The difference between applying this multilevel idea to least squares problems and applying it to stochastic path simulation is explained by the difference in trends of the variances. 
In applications with SDEs, the variance of $P_L$ does not vary with level, therefore the total cost of the simple Monte Carlo framework increases with level, for example $V_L^{MC}=\mathcal{O}(1), C_L^{MC} = \mathcal{O}(2^{L})$ so Monte Carlo achieves a total variance of $\varepsilon^2$ with cost $\mathcal{O}(\varepsilon^{-2} 2^{L})$.

Conversely, in our least squares application the variance of $x^{(L)}$ decreases with level as $\mathcal{O}(2^{-L})$, and the cost per sample at level $L$ is $\mathcal{O}(2^L)$, therefore simple Monte Carlo achieves total variance $\varepsilon^2$ with a total cost of $\mathcal{O}(\varepsilon^{-2})$. That is why the total cost does not depend on the level. The fact that $\varepsilon^2 \mathcal{C}_{MC}=\mathcal{O}(1)$ is independent of the level $L$ and only depends on the desired accuracy $\varepsilon$ is theoretically justified for Gaussian sketching mentioned in Section \ref{sec:varSimpleMC}, and we verified it numerically for the sketching methods mentioned in Section \ref{sec:BackgroundSketching}. 

In other words, in stochastic path simulation the multilevel approach allows us to reduce the overall cost of the estimation from $\mathcal{O}(\varepsilon^{-3})$ to  $\mathcal{O}(\varepsilon^{-2})$, whereas in least squares problems the overall cost is already $\mathcal{O}(\varepsilon^{-2})$ in simple Monte Carlo estimation, therefore a multilevel approach offers no savings.

\section{Discussion} 

We studied a multilevel approach to estimating the solution of large overdetermined least squares problems using random sketching or subsampling to obtain the individual random samples. In our estimator we double the sketch size each time we go up by one level. The analysis of the error shows a considerable reduction in the variance on each level, especially when antithetics is used. However the cost of generating samples using the QR factorization does not increase fast enough to give an advantage to our multilevel method in terms of total computational cost. 
Despite this, the multilevel approach may be of potential interest to achieve higher levels of parallelization, however we believe this would depend mainly on the implementation of the framework and is outside of the scope of this paper.

\subsection*{Acknowledgement}
The authors acknowledge the use of ChatGPT (OpenAI) to assist with proofreading and correction of spelling and grammar in the manuscript.

\bibliographystyle{plain}
\bibliography{refs}

@STRING(simax="SIAM J. Matrix Anal. Appl.")

@STRING(annals="Ann. Math.")

@article{nakatsukasa_tropp_2024fast,
  title={Fast and accurate randomized algorithms for linear systems and eigenvalue problems},
  author={Nakatsukasa, Yuji and Tropp, Joel A},
  journal=simax,
  volume={45},
  number={2},
  pages={1183--1214},
  year={2024},
  publisher={SIAM}
}

@article{murray2023randomized,
  title={Randomized numerical linear algebra: A perspective on the field with an eye to software},
  author={Murray, Riley and Demmel, James and Mahoney, Michael W and Erichson, N Benjamin and Melnichenko, Maksim and Malik, Osman Asif and Grigori, Laura and Luszczek, Piotr and Derezi{\'n}ski, Micha{\l} and Lopes, Miles E and others},
  journal={arXiv preprint arXiv:2302.11474},
  year={2023}
}

@misc{epperly_blogBias,
  author       = {Ethan N. Epperly},
  title        = {Note to Self: Sketch-and-Solve with a {Gaussian} Embedding},
  howpublished = {\url{https://www.ethanepperly.com/index.php/2024/11/19/note-to-self-sketch-and-solve-with-a-gaussian-embedding/}},
  month        = nov,
  year         = {2024},
  note         = {Accessed: 2025-08-21}
}

@article{bartan,
  title={Distributed sketching methods for privacy preserving regression},
  author={Bartan, Burak and Pilanci, Mert},
  journal={arXiv preprint arXiv:2002.06538},
  year={2020}
}

@article{sketched_ridgeReg,
author = {Wang, Shusen and Gittens, Alex and Mahoney, Michael W.},
title = {Sketched ridge regression: optimization perspective, statistical perspective, and model averaging},
year = {2017},
issue_date = {January 2017},
publisher = {JMLR.org},
volume = {18},
number = {1},
issn = {1532-4435},
journal = {J. Mach. Learn. Res.},
month = jan,
pages = {8039–8088},
numpages = {50}
}

@misc{precond_viaDistrib,
      title={Distributed Least Squares in Small Space via Sketching and Bias Reduction}, 
      author={Sachin Garg and Kevin Tan and Michał Dereziński},
      year={2024},
      eprint={2405.05343},
      archivePrefix={arXiv},
      primaryClass={cs.DS},
}

@article{giles2015actanum,
  title={Multilevel {Monte Carlo} methods},
  author={Giles, Michael B},
  journal={Acta Numerica},
  volume={24},
  pages={259--328},
  year={2015},
  publisher={Cambridge University Press}
}

@article{giles2008opre,
  title={Multilevel {Monte Carlo} path simulation},
  author={Giles, Michael B},
  journal={Operations research},
  volume={56},
  number={3},
  pages={607--617},
  year={2008},
  publisher={INFORMS}
}

@inproceedings{drineas2006sampling,
author = {Drineas, Petros and Mahoney, Michael W. and Muthukrishnan, S.},
title = {Sampling algorithms for $l_2$ regression and applications},
year = {2006},
isbn = {0898716055},
publisher = {Society for Industrial and Applied Mathematics},
address = {USA},
booktitle = {Proceedings of the Seventeenth Annual ACM-SIAM Symposium on Discrete Algorithm},
pages = {1127–1136},
numpages = {10},
location = {Miami, Florida},
series = {SODA '06}
}

@article{drineas2011faster,
  title={Faster least squares approximation},
  author={Drineas, Petros and Mahoney, Michael W and Muthukrishnan, Shan and Sarl{\'o}s, Tam{\'a}s},
  journal={Numerische Mathematik},
  volume={117},
  number={2},
  pages={219--249},
  year={2011},
  publisher={Springer}
}

@article{mahoney2011randomized,
  title={Randomized algorithms for matrices and data},
  author={Mahoney, Michael W},
  journal={Foundations and Trends{\textregistered} in Machine Learning},
  volume={3},
  number={2},
  pages={123--224},
  year={2011},
  publisher={Emerald Publishing Limited}
}

@article{woodruff2014sketching,
  title={Sketching as a tool for numerical linear algebra},
  author={Woodruff, David P},
  journal={arXiv preprint arXiv:1411.4357},
  year={2014}
}

@article{drineasRandNLA,
author = {Drineas, Petros and Mahoney, Michael W.},
title = {{RandNLA}: randomized numerical linear algebra},
year = {2016},
issue_date = {June 2016},
publisher = {Association for Computing Machinery},
address = {New York, NY, USA},
volume = {59},
number = {6},
issn = {0001-0782},
url = {https://doi.org/10.1145/2842602},
doi = {10.1145/2842602},
journal = {Commun. ACM},
month = may,
pages = {80–90},
numpages = {11}
}

@INPROCEEDINGS{sarlos,
  author={Sarlos, Tamas},
  booktitle={2006 47th Annual IEEE Symposium on Foundations of Computer Science (FOCS'06)}, 
  title={Improved Approximation Algorithms for Large Matrices via Random Projections}, 
  year={2006},
  volume={},
  number={},
  pages={143-152},
  doi={10.1109/FOCS.2006.37}}

@article{ma15a,
  author  = {Ping Ma and Michael W. Mahoney and Bin Yu},
  title   = {A Statistical Perspective on Algorithmic Leveraging},
  journal = {Journal of Machine Learning Research},
  year    = {2015},
  volume  = {16},
  number  = {27},
  pages   = {861--911},
  url     = {http://jmlr.org/papers/v16/ma15a.html}
}

@article{raskutti,
  author  = {Garvesh Raskutti and Michael W. Mahoney},
  title   = {A Statistical Perspective on Randomized Sketching for Ordinary Least-Squares},
  journal = {Journal of Machine Learning Research},
  year    = {2016},
  volume  = {17},
  number  = {213},
  pages   = {1--31},
  url     = {http://jmlr.org/papers/v17/15-440.html}
}

@article{pilanci2016iterative,
  title={Iterative Hessian sketch: Fast and accurate solution approximation for constrained least-squares},
  author={Pilanci, Mert and Wainwright, Martin J},
  journal={Journal of Machine Learning Research},
  volume={17},
  number={53},
  pages={1--38},
  year={2016}
}

@InProceedings{wang17d,
  title = 	 {{Sketching Meets Random Projection in the Dual: A Provable Recovery Algorithm for Big and High-dimensional Data}},
  author = 	 {Wang, Jialei and Lee, Jason and Mahdavi, Mehrdad and Kolar, Mladen and Srebro, Nati},
  booktitle = 	 {Proceedings of the 20th International Conference on Artificial Intelligence and Statistics},
  pages = 	 {1150--1158},
  year = 	 {2017},
  editor = 	 {Singh, Aarti and Zhu, Jerry},
  volume = 	 {54},
  series = 	 {Proceedings of Machine Learning Research},
  month = 	 {20--22 Apr},
  publisher =    {PMLR},
  url = 	 {https://proceedings.mlr.press/v54/wang17d.html}
}

@article{bartan2023distributed,
  title={Distributed sketching for randomized optimization: Exact characterization, concentration, and lower bounds},
  author={Bartan, Burak and Pilanci, Mert},
  journal={IEEE Transactions on Information Theory},
  volume={69},
  number={6},
  pages={3850--3879},
  year={2023},
  publisher={IEEE}
}

@article{acebron2020probabilistic,
  title={{A probabilistic linear solver based on a multilevel {Monte Carlo} method}},
  author={Acebr{\'o}n, Juan A},
  journal={Journal of Scientific Computing},
  volume={82},
  number={3},
  pages={65},
  year={2020},
  publisher={Springer}
}

@article{giles2014antithetic,
  title={Antithetic multilevel {Monte Carlo} estimation for multi-dimensional SDEs without L{\'e}vy area simulation},
  author={Giles, Michael B and Szpruch, Lukasz},
  year={2014},
  journal={The Annals of Applied Probability},
volume={24},
number={4},
pages={1585--1620}
}

@article{drineas2012levSco,
  title={Fast approximation of matrix coherence and statistical leverage},
  author={Drineas, Petros and Magdon-Ismail, Malik and Mahoney, Michael W and Woodruff, David P},
  journal={The Journal of Machine Learning Research},
  volume={13},
  number={1},
  pages={3475--3506},
  year={2012},
  publisher={JMLR. org}
}

@article{jones2025subapsnap,
  title={SubApSnap: Solving parameter-dependent linear systems with a snapshot and subsampling},
  author={Jones, Eleanor and Nakatsukasa, Yuji},
  journal={arXiv preprint arXiv:2510.04825},
  year={2025}
}

@article{tropp2011improved,
  title={Improved analysis of the subsampled randomized {H}adamard transform},
  author={Tropp, Joel A},
  journal={Advances in Adaptive Data Analysis},
  volume={3},
  number={01n02},
  pages={115--126},
  year={2011},
  publisher={World Scientific}
}

@book{golub2013matrix,
  title={Matrix computations},
  author={Golub, Gene H and Van Loan, Charles F},
  year={2013},
  publisher={JHU press}
}

\end{document}